\documentclass[12pt]{amsart}

\usepackage{graphicx}

\newtheorem{theorem}{Theorem}
\newtheorem{proposition}{Proposition}
\newtheorem{definition}{Definition}
\newtheorem{lemma}{Lemma}
\newtheorem{corollary}{Corollary}

\newcommand{\lcr}{\raisebox{-5pt}{\mbox{}\hspace{1pt}
                  \includegraphics{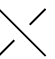}\hspace{1pt}\mbox{}}}
\newcommand{\ift}{\raisebox{-5pt}{\mbox{}\hspace{1pt}
                  \includegraphics{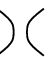}\hspace{1pt}\mbox{}}}
\newcommand{\zer}{\raisebox{-5pt}{\mbox{}\hspace{1pt}
                  \includegraphics{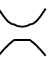}\hspace{1pt}\mbox{}}}

\title{Quantum Obstruction Theory}
\author{Charles Frohman}
\address{Department of Mathematics, University of Iowa, Iowa City, IA
52242, USA}
\email{\tt frohman@math.uiowa.edu}
\author{Joanna Kania-Bartoszy\'{n}ska}
\address{Department of Mathematics, Boise State University, Boise, ID
83725, USA}
\email{\tt kania@math.idbsu.edu}
\thanks{During the writing of this paper the first author was partially
  supported by NSF-DMS-9803233, and the second author by NSF-DMS-9626818.}

\begin{document}

\begin{abstract}
We use topological quantum field theory to derive an invariant of a
three-manifold with boundary.
We then show how to use the structure of this invariant as an obstruction
to embedding one three-manifold into another.
\end{abstract}

\maketitle

\section{Introduction}

In the mid 1980's the Jones polynomial of a link  was introduced \cite{jones1}, \cite{jones2}.
It was defined as the normalized trace of an element of a braid group, corresponding to the 
link, in a certain representation. 
Although a topological invariant, the extrinsic nature of its computation
made the relationship between topological configurations lying in the 
complement of the knot and the value of the Jones polynomial obscure.
In order to elucidate this connection, Witten \cite{Wi}
introduced Topological Quantum Field Theory.
TQFT is  a cut and paste technique which allows for localized computation of
the Jones polynomial.
Following his discovery various approaches to topological quantum
field theory were introduced (\cite{B}, \cite{FK}, \cite{G}, \cite{Law}, \cite{MS},
\cite{R}, \cite{Tu}, \cite{W} and others). In the work of Witten, it
is obvious that the vector 
spaces of topological quantum field theory are quantizations of the
space of connections on the manifold.
Initial considerations of TQFT focused on the formal rules of
combination of vector spaces and vectors.
As more intrinsic developments of
TQFT appeared there was a shift towards the use of the Kauffman bracket skein
module. The Kauffman bracket skein module of a 3-manifold can be
thought of as a quantization of the $SL_2({\mathbb C})$--characters of 
the manifold \cite {BFK}.
This approach to TQFT recovers some of the initial feel of Witten's work.

Although the study of $3$-manifold invariants is substantial in itself,
few applications to classical $3$-manifold topology have been found.
In this paper, we use quantum invariants
to develop obstructions
to embedding one 3-manifold in another. 
In order to
extract quantum invariants from a 3-manifold  a good deal
of extra information about its boundary  is needed. The big challenge
is to sort out what is information about the manifold and what is
just data about the parameterization chosen for the boundary. A very
similar problem arises in differential geometry. The Riemannian curvature tensor is a
classifying invariant of the metric structure of the manifold, but
it comes in the form of a tensor that depends intimately on the choice of
local coordinates. To sidestep this, invariant functions of the coefficients
of the tensor are used in order to obtain numbers that are independent
of the choice of local coordinates. 
Similarly, we take the quantum invariant of a manifold and perform
a construction to produce
an invariant that is independent of the extra information about how the boundary
is parametrized.

In the next section needed definitions are recalled  and the TQFT is constructed.  The path we
take is the elegant approach of Roberts \cite{R}, along with a simplification of
Sikora \cite{S}. This is a compact and eminently computable approach to
topological quantum field theory. The vector spaces are
built out of framed links in the manifold, so that they have a topological
feel much like singular homology and fundamental group.

In the following section 
we prove a theorem that gives an obstruction to embedding one 3-manifold
in another and illustrate it with an application of quantum
obstructions.  
If a 3-manifold contains a punctured lens space  then
it cannot be embedded in the 3-sphere. This observation lies at the heart of
the technique of Scharlemann cycles - a technique that has been extremely fruitful
in resolving questions about Dehn surgery on knots. Our proposition
generalizes this 
criterion and opens the door to the use of quantum invariants in similar combinatorial
arguments.

\section{Projective TQFT}

There is a construction of topological quantum field theory based on the
Witten-Reshetikhin-Turaev invariants due to Blanchet, Habegger, Masbaum and
Vogel \cite{B}. The vector space associated to a surface is a quotient of the vector
space with basis given by all connected oriented 3-manifolds having that surface as boundary. 
There is a pairing between
the vector space of a surface and the vector space associated to the
same surface 
with the opposite orientation.
Specifically, if $\partial M=-\partial N$ then $<M,N>=Z_r(M \cup N)$ where $Z_r$
is the Witten-Reshetikhin-Turaev invariant of level $r$ \cite{RT},
normalized so that the invariant 
of the 3-sphere is 1. To obtain the vector space
of a surface take the quotient of  the vector space above by the
radical of this form.
The invariant of a connected 3-manifold with boundary  is just its image in
the quotient.

We take a different approach here, with equivalent results, that makes
fundamental computations more direct.
The full apparatus of topological quantum field theory is not needed,
as the invariant we derive is independent of framing. Hence, we set up a 
projective topological quantum field theory. The approach we take is based
on chapter 6 of \cite{R}, clarified by using a theorem from
 \cite{S}. The construction uses the Kauffman bracket skein module
at a root of unity, modulo fusion.

\subsection{The Kauffman Bracket Skein Module}Let $r$ be an odd prime.
Let $u$ be a primitive 8rth root of unity and  $A=u^2$.  
We will work over 
$\mathbb{Q}[u]$, the cyclotomic numbers corresponding to $u$. 
The quantized integer  $[n]$  denotes 
\[ \frac {A^{2n}-A^{-2n}}{A^2-A^{-2}}, \]
which is  a unit in the ring $\mathbb{Z}[u]$ for $n<r$.
Let
$X$ be the positive square root of $\sum_{c=0}^{r-2}[c+1]^2$. It can
be shown that $X \in \mathbb{Z}[u]$.
Explicitly,
\[ X= \sqrt{\frac{r}{2}}\frac{1}{\sin{\frac{\pi}{r}}}.\]

Let $M$ be an orientable $3$-manifold.
A framed link in $M$ is an embedding of a disjoint union of annuli
into $M$. In our diagrams we  draw only the core of an annulus lying
parallel to the plane of the paper, i.e. with blackboard framing.

Two framed links in $M$ are equivalent if there is an isotopy of $M$
taking one to the other. Let $\mathcal{L}$ denote the set of equivalence
classes of framed links in $M$, including the empty link.  In the vector space  
$\mathbb{Q}[u] \mathcal{L}$, with  basis $\mathcal{L}$,
define $S(M)$ to be the smallest subspace of $\mathbb{Q}[u] \mathcal{L}$
containing all expressions of the form
$\displaystyle{\lcr-A\zer-A^{-1}\ift}$
and 
$\bigcirc+A^2+A^{-2}$,
where the framed links in each expression are identical outside the
balls pictured in the diagrams. The Kauffman bracket skein module
$K_r(M)$ is defined to be
the quotient
\[ \mathbb{Q}[u] \mathcal{L} / S(M). \]

We assume that the reader is familiar with
the Jones-Wenzl idempotents and the Kauffman triads as described on pages 136-138 of
 \cite{Li}. We  refer the reader to pages 150-152 of  \cite{Li} for a
 computation of the number $\theta(a,b,c)$ 
that is defined whenever $(a,b,c)$ is an admissible triple. Notice 
that  $\theta(a,b,c)$ is a unit in 
$ \mathbb{Z}[u]$.
By labeling an edge of a trivalent graph with the letter $m$ we mean 
that it carries the $m$-th Jones-Wenzl idempotent.

Fusion (see \cite{Li}, page 156) is the relation:
\[ \raisebox{-22pt}{\includegraphics{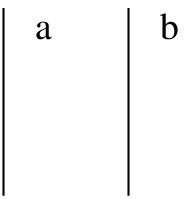}} \hspace{2pt}= 
\sum_{c} (-1)^c \frac{[c+1]}{\theta(a,b,c)} 
\hspace{2pt} \raisebox{-22pt}{\includegraphics{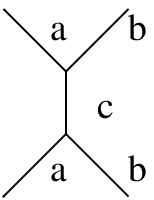}}, \]
where the sum is over all c so that the triads are admissible. The result of enforcing fusion
in the Kauffman bracket skein module $K_r(M)$ is denoted $K_{r,f}(M)$, and will
be referred to as the reduced Kauffman bracket skein module of level r.
The element $\Omega_r$ is the skein in the solid torus, 
\[ \frac{1}{X}\sum_{c=0}^{r-2} (-1)^c[c+1] \hspace{2pt} \raisebox{-15pt}
{\includegraphics{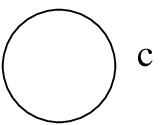}}.\]
When one plugs $\Omega_r$ into a component of a framed link in $M$, then in $K_{r,f}(M)$
one can alter the framed link by handle slides along that component and not change
the element of $K_{r,f}(M)$. The result of plugging  
in copies of  $\Omega_r$ into each of the components of the link $L$
is denoted by $\Omega_r(L)$.  Notice that the  Kauffman bracket of the
result of plugging  $\Omega_r$ into a trivial
framed knot in $S^3$ is $X$.

Let $\kappa$ be the Kauffman bracket of the skein
in the three sphere obtained by plugging $\Omega_r$ into the framed knot shown below.

\hspace{1.75in}\includegraphics{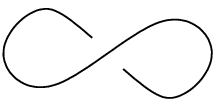}

We will also need the Turaev-Wenzl identity (see \cite{Li}). Suppose a 
diagram contains the figure below, where $\Omega_r$ is plugged into the horizontal circle,
and the vertical arc carries the $c$th Jones-Wenzl idempotent. The
diagram represents zero in the reduced Kauffman bracket skein module unless
$c=0$.

\begin{picture}(56,56)
\hspace{2in} \includegraphics{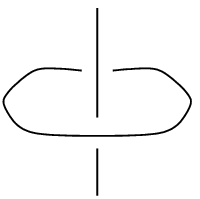}
\put(4,24){$\Omega_r$}
\put(-22,45){$c$}
\end{picture}

\subsection{preTQFT} We review a result of
Sikora that is used in the construction of the TQFT.  Suppose that $M$
and $M'$ are 
compact oriented 3-manifolds and $f: \partial M \rightarrow \partial
M'$
 is an orientation preserving
homeomorphism. We say that framed links $L \subset M$ and $L' \subset M'$ are adapted to $f$
if 
\begin{enumerate}
\item The result of surgery along $L$ is $M'$, the result of surgery along $L'$ is $M$.
\item There are regular neighborhoods  $N(L)$ and $N(L')$, and an extension of $f$, 
\[f:M-N(L) \rightarrow M'-N(L')\]
that is a homeomorphism  taking the meridian of $L$ to the framing of $L'$, and
the framing of $L$ to the meridian of $L'$.
\end{enumerate}

Given a pair of links adapted to $f$ and the extension, we can define a map from $K_{r,f}(M)$
to $K_{r,f}(M')$, which we denote
\[K_r(f): K_{r,f}(M) \rightarrow K_{r,f}(M').\]

Given a framed link $L_1$ in $M$ push it off of $L$, then map it forward to $M'$ using the
extension of $f$ to get $f(L_1)$.   The skein in $M'$ induced by the union 
of $f(L_1)$ and 
$\Omega_r(L')$ is  $K_r(f)(L_1)$. Any two links adapted to $f$ differ by the Kirby-Roberts 
moves for surgery on a manifold with
boundary \cite{R1}. From this we see that the map $K_r(f)$ is well defined up to a power of $\kappa$.

\begin{theorem} [\cite{S}] The map $K_r(f): K_{r,f}(M) \rightarrow K_{r,f}(M')$ 
is an isomorphism.\qed\end{theorem}

Let $M \# N$ denote the connected sum of $M$ and $N$. The Kauffman bracket skein
modules of $M \# N$ and of the disjoint union of $M$ and $N$ are not
isomorphic.
However, there
is a natural isomorphism between $K_{r,f}(M \# N)$  and $K_{r,f}( M
\sqcup N)$. Let
$M'$ and $N'$ be the result of removing a ball from each of $M$ and $N$. There
is a map $\iota :  M' \sqcup N' \rightarrow M \# N$ coming from inclusion. As
the Kauffman bracket skein module of the result of removing a ball from a manifold
is isomorphic to the Kauffman bracket skein module of the manifold, this induces,
\[ \hat{\iota}:K_{r,f}(M \sqcup N) \rightarrow K_{r,f}(M \# N).\]

\begin{proposition} The map $\hat{\iota}$ is an isomorphism.\end{proposition}

\proof 
This follows immediately from the following lemma.
\begin{lemma}
If a $3$-manifold $M$ contains a  family of embedded spheres then  any
skein in $K_{r,f}(M)$
can be represented by a linear combination of links that miss these
spheres. 
\end{lemma}

Let $S^2$ be a sphere embedded in $M$ and
suppose that $L$ is a link in $M$ which intersects $S^2$ transversely
in one point. One can isotope $L$ so it is the same away from the
sphere and has two added  twists close to the sphere. This implies that
$<L>=z<L>$ for some complex number $z\neq 1$, thus $<L>=0$.

Take a representative of any skein in $K_{r,f}(M)$. One can assume
that it is a linear combination of  trivalent graphs  admissibly colored with
Jones-Wenzl idempotents \cite{Lick}.   Use fusion
repeatedly to make sure that the graphs do not have more than one edge 
intersecting each of the spheres embedded in $M$. Use the light bulb
trick as above to see that the label on that edge has to be zero.
\qed

Hence, we
can be fluid about connected sums and disjoint unions so long as we work in the reduced
skein module.

\subsection{TQFT}
Recall that a {\em projective Topological Quantum Field Theory} consists of a
functor and a partition function, defined for each level $r$.
The functor $V_r$ maps  the
category of {\em marked} surfaces and homeomorphisms to the category of
vector spaces over $\mathbb{Q}[u]$, 
and linear maps defined up to scalar multiplication by a power of
$\kappa$. 
The partition function $Z_r$ assigns to
every {\em marked} 
3-manifold $M$ a vector $Z_r(M)$ in the vector space corresponding to
its boundary, well 
defined up to scalar multiplication by a power of $\kappa$.
These have to satisfy the following axioms.
\begin{enumerate}
\item {\bf Dimension}:
The empty surface is mapped to the vector space $\mathbb{Q}[u]$.
\[\emptyset \mapsto \mathbb{Q}[u]\]
\item {\bf Disjoint union}: 
Disjoint union of surfaces is mapped to the tensor product of vector
spaces.
\[F_1\sqcup F_2 \mapsto V_r(F_1)\otimes V_r(F_2)\]
\item {\bf Duality}: 
Oppositely oriented marked surfaces are mapped to naturally dual
vector spaces. 
\[V_r(-F)=V_r(F)^* \]
This duality determines a pairing:
\[\langle\ ,\ \rangle :V_r(F)\otimes V_r(-F)\rightarrow \mathbb{Q}[u]\]
\item {\bf Gluing}: 
Gluing two $3$-manifolds along surfaces in their boundaries corresponds 
to contracting vectors.

Let $M_1$ and $M_2$ be marked  $3$-manifolds with $\partial M_1 =F\sqcup F_1$
and $\partial M_2 =-F\sqcup F_2$ where the markings agree on $F$. Then  
$Z_r(M_1)\in V_r(F)\otimes V_r(F_1)$ can  be written as 
$\sum v\otimes v_1$, and $Z_r(M_2)\in V_r(F)^*\otimes V_r(F_2)$ can be
written as $\sum \phi\otimes v_2$. If $M=M_1\cup_F M_2$ then
\[Z_r(M)=\sum \phi(v) v_1\otimes v_2.\]
\item {\bf Mapping cylinder}:
The value of the partition function on the mapping cylinder of a homeomorphism
corresponds to a linear map given by that homeomorphism.

Let $f:F_1\rightarrow F_2$ and let $M=F_1\times I\cup_f F_2$, where
$(x,1)\equiv f(x)$, so that $\partial M= -F_1\sqcup F_2$. Then

\[Z_r(M)=V_r(f)\in V_r(F_1)^*\otimes V_r(F_2)\]

\end{enumerate}

We will now construct a projective TQFT. Fix a level $r$.

{\em Vector Spaces.} Let $F$ be a closed oriented surface. A {\em
  marking} of $F$ is a  
choice of a handlebody
$H$ with $\partial H=F$, up to homeomorphisms that are the identity on
$F$. We assign to  
the marked surface $(F,H)$ the vector space $V_r(F,H)=K_{r,f}(H)$. 

{\em Linear Maps.} Let $f: F \rightarrow F'$ be an orientation
preserving homeomorphism of  
the marked surfaces
$(F,H)$ and $(F',H')$. Given links $L\subset H$ and $L'\subset H'$ 
adapted to $f$, we let 
$V_r(f)=K_r(f)$. Our morphism is defined up to multiplication by a
power of $\kappa$. 

Suppose now that a surface is not connected. In this case, a marking
for the surface is a choice 
of  a handlebody $H_i$ for each component. The vector space associated with
the surface is $\otimes_i K_{r,f}(H_i)$. We can realize this tensor
product by taking 
any combination of connected sums and disjoint unions of the $H_i$,
since $K_{r,f}(\#_{i}H_i)$ is canonically isomorphic to 
$\otimes_i K_{r,f}(H_i)$.  

{\em Duality.} Let $(F,H)$ be a marked surface. The marked surface
with the orientations 
of $F$ and $H$ reversed is denoted by $-(F,H)$. We can form the double 
$H\cup_F -H$, 
which is homeomorphic to $\#_g S^1 \times S^2$, where $g$ is the genus
of $F$. 
Skeins $x \in K_{r,f}(H)$ and $y \in K_{r,f}(-H)$ determine a skein in
$K_{r,f}(\#_g S^1 \times S^2)$ obtained by taking their disjoint
union. In order to pair 
$x$ and $y$, choose a family of spheres that cut $\#_g S^1 \times S^2$ down to
punctured balls, and represent the union of $x$ and $y$ by a linear combination of links
missing the spheres. Since the corresponding links live in punctured balls
we can take their Kauffman bracket to get a complex
number times the empty link, finally multiply the coefficient by $X^g$. The final
multiplication makes
the answer coincide with a number computed in the 3-sphere.
To be more precise, choose an  isomorphism between 
$K_{r,f}(\#_g S^1\times S^2)$ 
and $K_{r,f}(S^3)$ coming from a surgery diagram, and the answer will
differ from the 
one we have given by multiplication by a power of $\kappa$. 
We have defined a pairing,
\[ <\ , \ > : V_r(F,H) \otimes V_r(-F,-H) \rightarrow \mathbb{Q}[u].\]
This is a duality pairing. Since $K_{r,f}(H)=K_{r,f}(-H)$, we
can see the pairing as defined on $V_r(F,H) \otimes V_r(F,H)$.

Choose a trivalent spine for the handlebody $H$. A coloring $a$ of the
spine  is a choice 
of an integer between $0$ and $r-2$ for each edge. When
the numbers at each vertex are admissible, we can form a skein by
running that 
number of arcs between the vertices and putting the appropriate
Kauffman triad at each 
vertex. We denote that skein by $x_a$. The skeins $x_a$ where the
coloring gives an 
admissible triad at each vertex form a basis for $K_{r,f}(H)$ \cite{Lick}. 
This can  easily be verified by
computing $<x_a,x_b>$. If $a\neq b$ the answer is $0$, and
$<x_a,x_a>= X^g u(a)$, where $u(a)$ is a unit in $\mathbb{Z}[u]$.
The dual basis $x^a$ are the skeins with $<x_b,x^a>=\delta^b_{a}$ where
$\delta^b_{a}$ is Kronecker's delta. It is easy to see that $x^a=\frac{1}{u(a)X^g}x_a$.

Suppose now that 
$C=\#H_i$.  There is a standard
embedding of
$C$ in $\#_g S^1\times S^2$, where $g$ is the sum of the genera of the
components of $\partial C$. 
In the case that $C$ is a handlebody, the double of $C$ is
homeomorphic to $\#_g S^1\times S^2$, 
and this is the standard embedding of $C$. 
If $C$ is a connected sum of handlebodies, take
the connected sum of the doubles of each of the handlebodies. Choose the
balls along which the 
sum is taken to lie inside the handlebodies  $H_i$. It is clear that
this embedding 
is unique up to homeomorphism and the complement of $C$ is a disjoint
union of handlebodies. 
Choose a trivalent spine for each $H_i$. A basis for $K_{r,f}(C)$ is
given by all admissible colorings of the
disjoint union of the spines, which we denote by $\otimes x_{a_i}$ or
$x_a$, where $a_i$ is the restriction
of the coloring $a$ to the $i$th spine.

We can see our pairing quite nicely in the three-sphere. To this
end we say that $C$ is embedded in $S^3$ in the standard way if its
complement is a  
disjoint union of handlebodies, so that each handlebody is unknotted
and there is a 
family of disjoint balls in $S^3$, each containing exactly one of the
complementary 
handlebodies.  Such an embedding is unique up to topological
equivalence. To get 
from a standard embedding of $C$ in $S^3$ to a standard embedding in
$\#_g S^1\times S^2$ 
we can do surgery on a link that lies in the complementary
handlebodies. Specifically, 
choose an unknotted circle with zero framing, linking each 1-handle of
$C$ lying in the 
complementary handlebodies so that the aggregate of such circles
forms the unlink which 
we denote $U$. These are the same as the circles used in the constructions
of \cite{MR}, \cite{R}. 
Surgery along this link yields  $\#_g S^1\times S^2$ with $C$ embedded
in a standard way.

In the three-sphere picture we can see the skeins $\otimes x_{a_i}$ lying in
$C$. To see a skein on the right hand side of the pairing, push it
through the boundary of $C$  and take its disjoint union  with
$\Omega_r(U)$. 
Notice that this is just the map of skein modules associated to the identity
on the boundary of  
the complementary handlebodies in  $\#_g S^1\times S^2$  to the
boundary of the  complementary handlebodies in $S^3$.

To pair $\otimes x_{a_i}$ with $\otimes x_{b_i}$, take the Kauffman bracket
of the union  $x_a\cup x_b \cup \Omega_r(U)$.  The
result of the pairing is
$\delta_{a}^b X^g  u(a)$, where $u(a)$ is a unit 
in $\mathbb{Z}[u]$ depending on the coloring $a$.

\subsection{Partition Function} Let $M$ be a compact, connected, oriented three-manifold. 
A marking for $M$ is a choice of
handlebodies, one for each connected component of $\partial M$. We denote
the connected sum of these handlebodies by $C$.
If $\partial M$ is empty then the marking is $S^3$. Since there is
only one choice we suppress it from the notation.
The {\em preinvariant}
of $(M,C)$ is the pair  consisting of $K_{r,f}(M)$ and the empty skein.
Let $L \subset M$ and $L' \subset C$ be framed links adapted to the identity
map on $\partial M$. Let $K_r(Id) : K_{r,f}(M) \rightarrow K_{r,f}(C)$ be the associated
isomorphism. We define 
\[Z_r(M,C)=K_r(Id)(\emptyset).\]
Once again, from the Kirby calculus
for manifolds with boundary,
this is well defined up to multiplication by a power of $\kappa$. Notice that the
orientation on $\partial M$ is always the orientation inherited from $M$.
Using the 
main result from \cite{MR} we see that if $M$ is any manifold then $Z_r(M,C)$ can be written as
\[ \frac{1}{X^g} \sum_a I_a \otimes x_{a_i},\] 
where the $I_a$ are algebraic integers.  The dual basis allows us to
rewrite this answer without any fractions:

\[Z_r(M,C)=\sum_a J_a \otimes x^{a_i},\]
where the $J_a$ are just the $I_a$ times a unit depending on $a$.

If $\partial M=\emptyset$ then $C=S^3$ and the basis consists of the
empty skein. In this case $Z_r(M)\in \mathbb{Q}[u]$ and, up to
multiplication by a power of $\kappa$, it is equal to the
Witten-Reshetikhin-Turaev invariant $\tau_r(M)$ \cite{KM}, \cite{RT}.

\subsection{Gluing.} Suppose that $(M_1,C_1)$ and $(M_2,C_2)$ are
marked $3$-manifolds, and
that there is a surface $F$ with $F\subset\partial M_1$ and
$-F\subset\partial M_2$. We will assume that the genus of $F$ is $g$. Suppose further 
that the handlebodies in $C_1$ and $C_2$
corresponding to $F$ and $-F$ are a handlebody $H$ for $F$ and $-H$ for $-F$.
The invariants of the manifolds are 
\[  \frac{1}{X^{g}} \sum_a I_a x_{a_1} \otimes \left(\otimes_{i>1} x^{a_i}\right),\]
and 
\[   \sum_b J_b x^{b_1} \otimes \left(\otimes_{j>1} x^{b_j}\right),\]
where  the first factor in each tensor corresponds to $F$ or $-F$. The invariant
of the result of gluing $(M_1,C_1)$ and $(M_2,C_2)$ with respect to the marking
$C$, coming from deleting the handlebodies $H$ and $-H$ and taking the
disjoint union
of the remaining handlebodies,
is 
\[  \frac{1}{X^{g}} \sum_{a,b} I_a J_b 
<x_{a_1},x^{b_1}>\left(\otimes_{i>1} x^{a_i}\right) \otimes\left( \otimes_{j>1} x^{b_j}\right).\]

To see this, glue $C_1$ to $C_2$ and then choose the obvious surgery diagram in
$-H$ to take the union to a connect sum of handlebodies corresponding to the
deletion of $H$ and $-H$. 

Using ideas of  Roberts \cite{R}, we can establish the equivalence of this projective
TQFT and the one derived from the constructions in \cite{B}. There is a natural
equivalence of the functors from the compact marked oriented surfaces to vector
spaces, that takes every manifold $M$ to a surgery diagram for $M$
in the handlebody that is its marking.
It is easy to check that this equivalence preserves the invariants of the three-manifolds
(projectively).

\section{Ideals and Obstructions}

In this section we introduce three ideals in $\mathbb{Z}[u]$  derived
from quantum invariants. The first two are invariants of the manifold, 
and can be used as obstructions to embedding one manifold into another.
The third ideal is not a manifold invariant but it can be easily
computed and is used to estimate one of the first two ideals.
If an ideal is all of $\mathbb{Z}[u]$ we say it is trivial. It is clear
that an ideal is trivial if it contains a unit \cite{ma}.

Let $M$ be a compact, oriented $3$-manifold with boundary. 

\begin{definition}
\[I_r(M)=\langle Z_r(M\cup N)\rangle_{\partial M =-\partial N}
 \leq {\mathbb Z}[u]\]
That is, $I_r(M)$ is the ideal of ${\mathbb Z}[u]$  generated by the
Witten-Reshetikhin-Turaev invariants of all closed manifolds
containing $M$.
\end{definition}
Similarly,
\begin{definition}
Let
\[I^{TV}_r(M)=
\langle Z_r\left((M\cup N)\#(-M\cup-N)\right)\rangle_{\partial M =-\partial N}
 \leq {\mathbb Z}[u]\]
be the ideal generated by the  Turaev-Viro invariants of all closed
manifolds containing $M$.
\end{definition}

\begin{proposition} The ideal $I_r(M)$ is an invariant of oriented homeomorphism.
The ideal $I^{TV}_r(M)$ is a topological invariant.
\qed
\end{proposition}

These are related in the following way.
\begin{proposition}
\begin{equation}
I_r^{TV}(M)\subset I_r(M)
\end{equation}
\end{proposition}
\proof
Any element of $I_r^{TV}(M)$ can be written as a sum, taken over all
the manifolds $N$ with $\partial N = -\partial M$,
\[\sum_N \alpha_N Z_r(M\cup N) \overline{Z_r (M\cup N)}.\] 
Since for every such $N$, the element $Z_r(M\cup N)$ is in the ideal $I_r(M)$,
thus also $Z_r(M\cup N) \overline{Z_r (M\cup N)}\in I_r(M)$, and so is the linear 
combination of such elements.
\qed

In particular, if $I_r(M)$ is nontrivial then $I_r^{TV}(M)$ is also
nontrivial.

\begin{proposition}\label{obstruct}
If the compact, oriented 3-manifold $M$ embeds in the compact, oriented 3-manifold $N$,
then $I_r(N) \subset I_r(M)$ and  $I^{TV}_r(N) \subset I^{TV}_r(M)$.

\end{proposition}

\proof
\[
I_r(N)=\langle Z_r(N\cup N')\rangle_{\partial N'=-\partial N} =
\langle Z_r\left((M\cup \text{Cl}(N-M))\cup N'\right)
\rangle_{\partial N'=-\partial N}
\]

\[
\subset \langle Z_r(M\cup N'')\rangle_{\partial N''=-\partial M}
\]

The argument for $I^{TV}_r$  is analogous.
\qed

\begin{corollary}\label{summand} 
Suppose that $N$ is a compact oriented $3$-manifold which contains a
submanifold $M$ such that $I_r(M)$ or $I_r^{TV}(M)$ is nontrivial,
then $N$ does not embed in $S^3$. \qed 
\end{corollary}

\begin{corollary}
If $M$ is a compact oriented $3$-manifold which contains a punctured lens space
$L(p,q)$, with $p\neq 2^n$ for $n>1$, then $M$ cannot be embedded into $S^3$.
\end{corollary}

\proof
It follows form
Theorem 2 of Yamada \cite{Y}, where he computes $Z_r(L(p,q))$. In
particular for any lens space there exists $r$ for which  $Z_r(L(p,q))$
is not a unit. If $p=2$  then for any odd prime $r$ the ideal is
zero. Otherwise take $r$ to be an odd prime dividing $p$, 
and the ideal will be either $0$ or nontrivial.
\qed

\begin{definition}
Let $(M,C)$ be a marked $3$-manifold, denote by
$g$ the sum of the genera of the handlebodies making up $C$, and suppose
that
\[Z_r(M,C)=\frac{1}{X^g}\sum_ai(a)x_a,\] with respect to a basis for $C$ made
up of admissibly colored trivalent spines. Define the ideal $J_r(M,C)$ 
to be generated by all coefficients $i(a)$.
\begin{equation}
J_r(M,C)= \langle i(a)\rangle \leq {\mathbb Z}[u].
\end{equation}
\end{definition}

The ideal
$J_r(M,C)$ is only an invariant of the marked three-manifold.
For instance, if we mark the solid torus $S^1 \times D^2$
with itself, then the surgery diagram is empty and the ideal is  generated
by $X$. If we mark the solid torus so that its longitude  corresponds
to the meridian of the marking solid torus and its meridian
corresponds to the longitude of the marking torus, 
then the surgery diagram consists of the core of the marking torus with $\Omega_r$ plugged
in, hence the ideal is 1.
We are
less interested in $J_r(M,C)$  as an invariant than we are as a  tool to estimate
$I^{TV}_r(M)$.

We restrict our attention now to the case when $\partial M$ is connected.
The ideal $I_r^{TV}(M)$ can be estimated from the ideal $J_r(M\#_D -M,C)$ for 
a particular marking $C$. Here $M\#_D -M$ denotes the disk sum of $M$
and $-M$ along a disk in $\partial M$. The construction is an extension 
for analyzing the invariant of a $3$-manifold with boundary of
Roberts' ``chain mail'' \cite{R}. 

\begin{proposition}\label{mark} If $\partial M$ is connected then
there is a marking $C$ for which
\[ I_r^{TV}(M) \subset \frac{1}{X^g} J_r(M\#_D -M,C).\]
\end{proposition}
\proof
Let $F=\partial M$ be a surface of genus $g$. Let
$D\subset F$ be a disk and let $F'=\overline{F-D}$, and denote the double
of $F'$ by
$DF=F' \cup_{\partial F'} -F'$. The boundary of $M\#_D -M$ can be naturally identified with 
$DF$. We will use the handlebody $F'\times I$ for the marking $C$ of
$M\#_D -M$, obviously $\partial (F'\times I)=DF$.

We will show that for any $N$ with $\partial M = - \partial N$, the Turaev-Viro invariant of
the manifold $M \cup_{\partial M=-\partial N} N$ can be expressed as an integral linear combination of the coefficients
$i(a)$ where,
\[ Z_r(M\#_{D} -M, C)= \frac{1}{X^g}\sum_{a} i(a)x_a.\]
Notice that $g$ is one half of the genus of $C$.
We achieve this by constructing an efficacious surgery description in
the three-sphere.
First we find a surgery diagram in $S^3$ for the double of $F' \times I$. Next, we locate a
framed link in $F' \times I$ that is adapted to the marking of $\partial (M\#_{D} -M)$
by $DF$. Finally we combine these diagrams to effect the computation.

The double of $F'\times I$ can be constructed in two steps.
First identify $F'\times \{1\}$ with $F'\times \{0\}$.
Second, fill the boundary so that the curves $\{pt\}\times I$, where
$pt\in \partial F'$, bound meridian disks in the filling torus.
Alternatively, one can do $0$-surgery along a fiber of $F\times S^1$.
A surgery diagram for $F\times S^1$ is drawn below for genus $g=2$. In
general, the number of clasped pairs woven through the middle circle
is equal to $g$.

\hspace{1.25in} \includegraphics{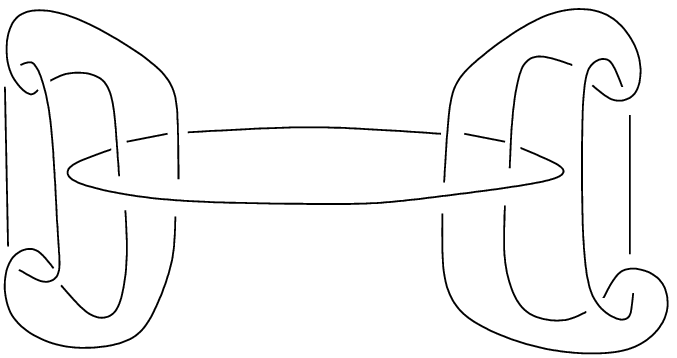}

Surgery along a fiber amounts to adding a zero framed unknotted circle that
links the middle circle.

\hspace{1.25in} \includegraphics{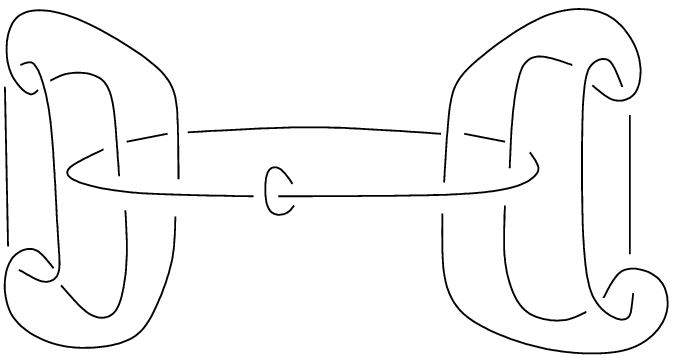}

The circumcision move says that if the small unknotted circle
and the middle circle are removed, then  surgery along the resulting diagram yields the
same result.

\hspace{1.25in} \includegraphics{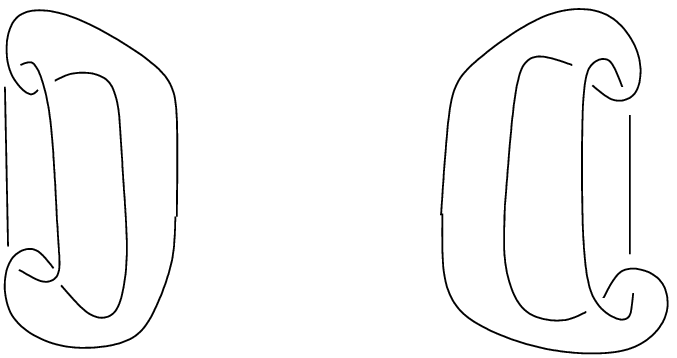}

Let $N(\partial D)$ be a regular neighborhood of the boundary of the
disk $D$ in $F$.
Notice that $N(\partial D) = S^1 \times [0,1]$.
Let $f: F\rightarrow {\mathbb R}$ be a smooth function, equal to $1$ on
the disk $D-N(\partial D)$, and $0$ outside of $D \cup N(\partial D)$. Suppose
that on $N(\partial D)$ the function $f$  has 
the circles $ S^1 \times \{pt\}$ as level sets, 
and that it decreases from $1$ to $0$ as you move from the boundary component
in $D$ to the boundary component outside of $D$, without any critical points.
Extend $f$ to a smooth function on $M$.
Perturb $f$ so that it is Morse and then move all index $1$ critical
points below all of the  index $2$ critical points, and cancel all local minima
and maxima. Let $l$ be a level between the index $1$ and index $2$
critical points. 
Notice that $f^{-1}([0,l])$ is the result of adding $1$ handles to
$F'\times I$ along $F'\times \{1\}$. We get $M$ from $f^{-1}([0,l])$
by adding $2$-handles along some curves $\gamma_i'$ in
$f^{-1}(l)$. Using this picture we can build a surgery diagram for $M\#_D-M$.
Let $G$ be the result of adding $h$ trivial $1$-handles from above to 
$F'\times \{\frac{1}{2}\}$, so that it has the same genus as $f^{-1}(l)$.

\begin{picture}(112,110)
\hspace{1.75in} \includegraphics{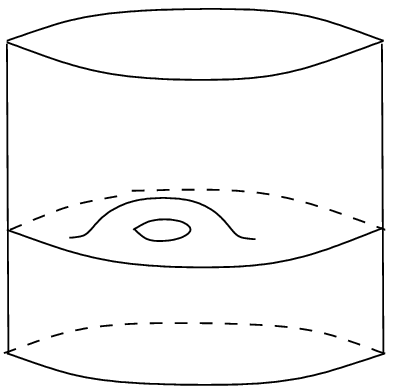}
\put(5,5){$F'$}
\put(5,45){$G$}
\put(5,95){$D$}
\end{picture}

Choose a homeomorphism from $f^{-1}([0,l])$ to the part of $F'\times I$ 
between $F'\times \{0\}$ and $G$, so that $f|_{F'\times \{0\}}$ is the
  identity.
Let $\gamma_i$ be the image of
the curves along which you add the $2$-handles. The surgery diagram will
consist of $g+h$ curves $\gamma_i$, along with $h$ $0$-framed curves
$c_{\lambda}$ above $G$ linking each added trivial $1$-handle in
$F'\times I$.

\begin{picture}(112,110)
\hspace{1.75in}\includegraphics{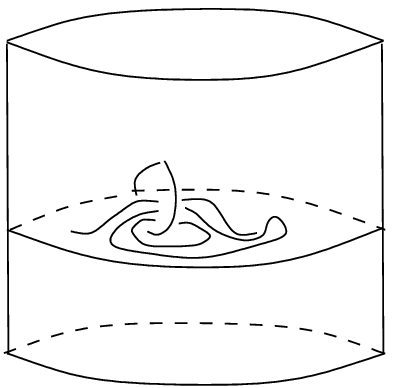}

\put(-30,45){$\gamma_i$}
\put(-60,60){$c_{\lambda}$}
\end{picture}

Notice that there is an orientation reversing map of order $2$ on the
result of surgery along the curves $c_{\lambda}$, that is the identity on $G$. This
establishes
that surgery on this diagram yields $M\#_D -M$. 

From this diagram and the marking $C$ we can compute
\begin{equation}\label{msum}
Z_r(M\#_D -M,C)=\frac{1}{X^{g}} \sum_a i(a)x_a,
\end{equation}
where  the $x_a$ are  standard basis elements coming from an admissibly
colored trivalent spine for $F'$ and $g$ is the genus of $F'$. To do
this, apply the Turaev-Wenzl identity $h$ times to cancel out the curves $c_{\lambda}$.

To see $G$ in the  surgery diagram for the double of $F' \times I$, take a disk 
whose boundary is where the
middle  curve of our initial surgery diagram was, and let it swallow the
arcs  of the clasped curves  going through it upwards. Finally add enough
trivial handles above to obtain a surface homeomorphic to $G$. 
Visualize the surgery diagram
for $M\#_D -M$ by seeing the curves $\gamma_i$ lying on $G$ and the
curves $c_{\lambda}$ linking the trivial handles.

\begin{picture}(191,66)
\hspace{1.25in} \includegraphics{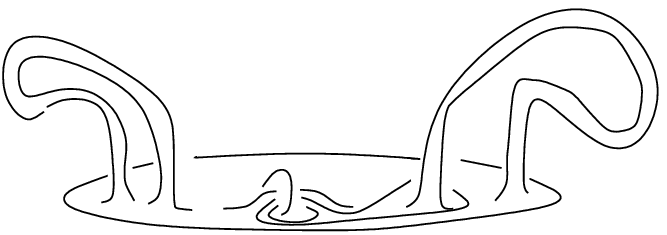}
\put(-60,14){$\gamma_i$}
\put(-105,17){$c_{\lambda}$}
\end{picture}

Since a regular neighborhood of the version of $G$ lying in $S^3$, along with balls containing
the trivial one handles and the linking curves, is homeomorphic to $F' \times I$,
we can expand the surgery curves $\gamma_i$ and the trivial linking handles
$c_{\lambda}$ to get the same formal sum $\frac{1}{X^{g}} \sum_a i(a)x_a$
as above, where the $x_a$ now lie in $S^3$.
 
We can perform an analogous construction based on a disk outside  the middle curve,
only swallowing the clasping curves downwards, and add trivial handles downwards.  If
$N$ is any $3$-manifold with $\partial N = -F$ we can find a surgery diagram for the
disk sum of $N$ with $-N$ from a Heegaard surface, and construct a surface $G'$ analogous
to $G$. We can place this surgery diagram on the surface outside that we just constructed
(with linking curves below the trivial handles). In the picture we just show the
tubes that swallow the clasping curves below. Trivial handles below can be added
as needed.

\hspace{.6in} \includegraphics{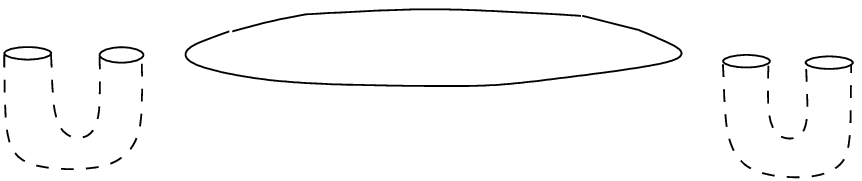}

In this way we can build a surgery diagram for any 
\[(M\#_D-M)\cup_{-\partial M = \partial N}(N\#_D-N).\]
Expand the
surgery curves in terms of a standard basis $x_a'$ associated to 
$F'\times I$, with $x'_a$ lying in the copy of $F'$ outside the middle 
curve.
Similar to equation (\ref{msum}),
\begin{equation}
Z_r(N\#_D -N,C)=\frac{1}{X^{g}} \sum_a j(a)x'_a.
\end{equation}
The $2g$ clasped  circles are trivial and unlinked.
Furthermore, the disks these circles bound
cut the Heegaard surfaces down to planar surfaces. Hence,  the Turaev-Wenzl
identity can be applied enough times to get rid of all the fractional part of the invariant.
Therefore,
\[Z_r\left((M\#_D-M)\cup (N\#_D-N) \right)\]
is an integral linear combination 
of $i(a)j(b)$. This implies that 
\[I_r^{TV}(M) \subset\frac{1}{X^g} J_r(M\#-M, C)\]
 where $C$
is the marking we constructed here.
\qed

As an immediate consequence of Propositions \ref{mark} and \ref{obstruct} we see that if 
the ideal $\frac{1}{X^g}J_r(M\#_D-M,C)$ is nontrivial then $M$ does not embed in the $3$-sphere. 
This leads to the following example.

Let $M$ be the result of gluing two solid tori along annuli in their boundaries.
Choose the annuli so that their cores are $(2,1)$ curves. The disk sum of $M$ 
with $-M$ can  be realized as surgery on a diagram in the cylinder over a punctured torus. 
The ideal   $\frac{1}{X}J_3(M\#_D-M,C)$ is the principal ideal generated by the integer $2$
 in the ring $Z[u]$. As this ideal is nontrivial, $M$ cannot be embedded in $S^3$.

\end{document}